\documentclass{article}

\usepackage{amssymb,latexsym,amsmath}

\usepackage{graphicx}

\begin{document}

\newcommand{\bfi}{\bfseries\itshape}

\makeatletter

\@addtoreset{figure}{section}

\def\thefigure{\thesection.\@arabic\c@figure}

\def\fps@figure{h, t}

\@addtoreset{table}{bsection}

\def\thetable{\thesection.\@arabic\c@table}

\def\fps@table{h, t}

\@addtoreset{equation}{section}

\def\theequation{\thesubsection.\arabic{equation}}

\makeatother

\newtheorem{theorem}{Theorem}[section]

\newtheorem{proposition}[theorem]{Proposition}

\newtheorem{lema}[theorem]{Lemma}

\newtheorem{cor}[theorem]{Corollary}

\newtheorem{definition}[theorem]{Definition}

\newtheorem{remark}[theorem]{Remark}

\newtheorem{exempl}{Example}[section]

\newenvironment{exemplu}{\begin{exempl}  \em}{\hfill $\square$

\end{exempl}}

\newcommand{\comment}[1]{\par\noindent{\raggedright\texttt{#1}

\par\marginpar{\textsc{Comment}}}}

\newcommand{\todo}[1]{\vspace{5 mm}\par \noindent \marginpar{\textsc{ToDo}}\framebox{\begin{minipage}[c]{0.95 \textwidth}

\tt #1 \end{minipage}}\vspace{5 mm}\par}

\newcommand{\ea}{\mbox{{\bf a}}}

\newcommand{\eu}{\mbox{{\bf u}}}

\newcommand{\ueu}{\underline{\eu}}

\newcommand{\ueo}{\overline{u}}

\newcommand{\oeu}{\overline{\eu}}

\newcommand{\ew}{\mbox{{\bf w}}}

\newcommand{\ef}{\mbox{{\bf f}}}

\newcommand{\eF}{\mbox{{\bf F}}}

\newcommand{\eC}{\mbox{{\bf C}}}

\newcommand{\en}{\mbox{{\bf n}}}

\newcommand{\eT}{\mbox{{\bf T}}}

\newcommand{\eL}{\mbox{{\bf L}}}

\newcommand{\eR}{\mbox{{\bf R}}}

\newcommand{\eV}{\mbox{{\bf V}}}

\newcommand{\eU}{\mbox{{\bf U}}}

\newcommand{\ev}{\mbox{{\bf v}}}

\newcommand{\eve}{\mbox{{\bf e}}}

\newcommand{\uev}{\underline{\ev}}

\newcommand{\eY}{\mbox{{\bf Y}}}

\newcommand{\eK}{\mbox{{\bf K}}}

\newcommand{\eP}{\mbox{{\bf P}}}

\newcommand{\eS}{\mbox{{\bf S}}}

\newcommand{\eJ}{\mbox{{\bf J}}}

\newcommand{\eB}{\mbox{{\bf B}}}

\newcommand{\eH}{\mbox{{\bf H}}}

\newcommand{\leb}{\mathcal{ L}^{n}}

\newcommand{\eI}{\mathcal{ I}}

\newcommand{\eE}{\mathcal{ E}}

\newcommand{\hen}{\mathcal{H}^{n-1}}

\newcommand{\eBV}{\mbox{{\bf BV}}}

\newcommand{\eA}{\mbox{{\bf A}}}

\newcommand{\eSBV}{\mbox{{\bf SBV}}}

\newcommand{\eBD}{\mbox{{\bf BD}}}

\newcommand{\eSBD}{\mbox{{\bf SBD}}}

\newcommand{\ecs}{\mbox{{\bf X}}}

\newcommand{\eg}{\mbox{{\bf g}}}

\newcommand{\paromega}{\partial \Omega}

\newcommand{\gau}{\Gamma_{u}}

\newcommand{\gaf}{\Gamma_{f}}

\newcommand{\sig}{{\bf \sigma}}

\newcommand{\gac}{\Gamma_{\mbox{{\bf c}}}}

\newcommand{\deu}{\dot{\eu}}

\newcommand{\dueu}{\underline{\deu}}

\newcommand{\dev}{\dot{\ev}}

\newcommand{\duev}{\underline{\dev}}

\newcommand{\weak}{\stackrel{w}{\approx}}

\newcommand{\mild}{\stackrel{m}{\approx}}

\newcommand{\lrightarrow}{\stackrel{L}{\rightarrow}}

\newcommand{\rrightarrow}{\stackrel{R}{\rightarrow}}

\newcommand{\strong}{\stackrel{s}{\approx}}

\newcommand{\weakdown}{\rightharpoondown}

\newcommand{\opg}{\stackrel{\mathfrak{g}}{\cdot}}

\newcommand{\opunu}{\stackrel{1}{\cdot}}
\newcommand{\opdoi}{\stackrel{2}{\cdot}}

\newcommand{\opn}{\stackrel{\mathfrak{n}}{\cdot}}
\newcommand{\opx}{\stackrel{x}{\cdot}}

\newcommand{\tr}{\ \mbox{tr}}

\newcommand{\Ad}{\ \mbox{Ad}}

\newcommand{\ad}{\ \mbox{ad}}

\renewcommand{\contentsname}{ }

\title{Braided spaces with dilations and sub-riemannian symmetric spaces}

\author{Marius Buliga \\
\\
Institute of Mathematics, Romanian Academy, 
P.O. BOX 1-764, \\
 RO 014700, Bucure\c sti, Romania\\
{\footnotesize Marius.Buliga@imar.ro}}

\date{(Final version, appeared in: Geometry. Exploratory Workshop on
Differential Geometry and its Applications, eds. D. Andrica, S. Moroianu,
Cluj-Napoca 2011, 21-35)}

\maketitle

\begin{abstract}
Braided sets which are also spaces with dilations 
are presented and explored in this paper, in the general frame of 
emergent algebras. Examples of such spaces are the sub-riemannian 
symmetric spaces.  
\end{abstract}

\vspace{.5cm}






\section{Introduction}

In the previous paper \cite{buligairq} we introduced and studied emergent
algebras, as a generalization of  differentiable algebras. An emergent algebra 
is a uniform idempotent right quasigroup, definition \ref{deftop}.

In this paper we explain with details previous results concerning 
conical groups, section \ref{secon}, dilatation structures 
(metric spaces with dilations), section \ref{dilats} as well as
new results concerning braided sets which are also dilatation structures,  
section \ref{secbraid}, in the frame of emergent algebras. Finally, we show that sub-riemannian
symmetric spaces (which are not Loos symmetric spaces) can be seen as braided 
dilatation structures. 

There is another, but related, line of research concerning symmetric spaces as 
emergent algebras, based on the notion of a "approximate symmetric space".
We postpone the presentation of this for a future paper. 

For yet another research line concerning the generalization of spaces with
dilations to deformations of normed groupoids see the paper 
\cite{buligagr}.

\paragraph{Acknowledgements.} Thanks are due to Radu Iordanescu for the
invitation 
to participate to the Workshop on Differential Geometry and its Applications, 
Iasi 2009, where part of this research was presented. I also wish to 
express my thanks to IHES, where another part of this work
has been done during a visit in March 2010. 

\section{Quandles, Loos symmetric spaces, contractible groups}

For braided sets see the paper \cite{ess}. 

\begin{definition}
Let $X$ be a non-empty set and $S: X \times X \rightarrow X \times X$ be a
bijection, 
$\displaystyle S(x_{1}, x_{2}) = (S_{1}(x_{1}, x_{2}), S_{2}(x_{1},x_{2}))$. We define for  $i=1,2$ the maps $\displaystyle S^{i i+1}: X^{3}
\rightarrow X^{3}$, $\displaystyle S^{12} = S \times \, id_{X}$, $S^{23} = \, 
id_{X} \times S$. 

(i) A map $S$ is called non-degenerate if for any fixed $y,z \in X$ the maps 
$x \mapsto S_{2}(x,y)$ and $x \mapsto S_{1}(z,x)$ are bijections.

(ii) A pair $(X,S)$ is a braided set (and $S$ is a braided map) if $S$ satisfies the braid relation
\begin{equation}
S^{12} S^{23} S^{12} \, = \, S^{23} S^{12} S^{23} 
\label{braideq}
\end{equation}

(iii) $S$ is involutive if $\displaystyle S^{2} = \, id_{X \times X}$. 

A braided set which is involutive is called symmetric set. 
\label{dess}
\end{definition}

A  large class of braided sets is made by pairs $(X,S)$, with 
$$S(x,y) = (x * y, x)$$
and $(X,*)$ is a rack. Racks and quandles are right quasigroups, a notion that
we shall use further, so here is the definition. 

\begin{definition} A right quasigroup is a set $X$ with a binary operation 
$*$ such that for each $a, b \in X$ there exists a unique $x \in X$ such that 
$a \, *  \, x \, = \, b$. We write the solution of this equation 
$x \, = \, a \, \backslash \, b$. 

A quasigroup is a set $X$ with a binary operation 
$*$ such that for each $a, b \in X$ there exist  unique elements  $x, y \in X$ 
such that $a \, *  \, x \, = \, b$ and $y \, * \, a \, = \, b$. 
We write the solution of the last  equation 
$y \, = \, b \, / \, a$. 

 An idempotent right quasigroup (irq) is a  right quasigroup $(X,*)$ such that 
 for any $x \in X$ $x \, * \, x \, = \, x$. Equivalently, it can be seen as a 
  set $X$ endowed with two  operations $\circ$ and $\bullet$, which satisfy the following axioms: for any $x , y \in X$  
\begin{enumerate}
\item[(P1)] \hspace{2.cm} $\displaystyle x \, \circ \, \left( x\, \bullet \,  y \right) \, = \, x \, \bullet \, \left( x\, \circ \,  y \right) \, = \, y$
\item[(P2)] \hspace{2.cm} $\displaystyle x \, \circ \, x \, = \, x \, \bullet \, x \,  = \,  x$
\end{enumerate}
The correspondence between notations using $*, \, \backslash$ and those using 
$\circ, \, \bullet$, is: $* \, = \, \circ$, $\backslash \, = \, \bullet$. 
\label{defquasigroup}
\end{definition}

In knot theory,  J.C. Conway and 
G.C. Wraith, in their unpublished correspondence  from 1959,  used the name 
"wrack" for  a  self-distributive right quasigroup generated by a link diagram. Later, Fenn and Rourke
\cite{fennrourke} proposed the name "rack" instead. Quandles are particular case
of racks, namely self-distributive idempotent 
right quasigroups. They were introduced by 
 Joyce \cite{joyce}, as a distillation  of the  Reidemeister
moves. More precisely, the axioms of a (rack ;  quandle ; irq) correspond respectively to the (2,3 ; 1,2,3 ; 1,2) Reidemeister moves.

We are interested in two particular cases of quandles. The first is related to 
symmetric spaces in the sense of Loos  \cite{loos} chapter II, definition 1. 


\begin{definition}
 $(X, inv)$ is a Loos algebraic symmetric space if $inv : X \times X \rightarrow X$ 
is an operation which satisfies the following axioms: 
\begin{enumerate}
\item[(L1)] $inv$ is idempotent: for any $x \in X$ we have 
$ inv (x,x) = \, x$, 
\item[(L2)] distributivity: for any $x, y , z \in X$ we have 
$$\displaystyle inv (x ,  inv \left( y ,  z \right)) \, = \, 
inv \left( inv \left( x,  y \right) , inv \left( x ,  z \right) \right)$$ 
\item[(L3)] for any $x, y \in X$ we have $inv \left( x , inv 
\left( x  ,  y \right) 
\right) \, = \, y$, 
\item[(L4)] for every $x \in X$ there is a neighbourhood $U(x)$ such that 
 $inv (x ,  y)  \, = \, y$  and $y \in U(x)$ then $x = y$. 
\end{enumerate}
If $X$ is a manifold, $inv$ is smooth (of class $\displaystyle 
\mathcal{C}^{\infty}$) and (L4) is true locally then $(X, inv)$ is a symmetric
space as defined by Loos \cite{loos} chapter II, definition 1. 
\label{defsym}
\end{definition}

Remark that if $(X,\, inv)$ is a Loos symmetric space then it is clearly a
quandle, therefore $(X,Inv)$ is a 
braided symmetric set, where: 
$$Inv(x,y) \, = \, (inv^{x} y , x)$$

\begin{definition}
A contractible group is a pair $(G,\alpha)$, where $G$ is a  
topological group with neutral element denoted by $e$, and $\alpha \in Aut(G)$ 
is an automorphism of $G$ such that: 
\begin{enumerate}
\item[-] $\alpha$ is continuous, with continuous inverse, 
\item[-] for any $x \in G$ we have the limit $\displaystyle 
\lim_{n \rightarrow \infty} \alpha^{n}(x) = e$. 
\end{enumerate}
\label{defunu}
\end{definition}

If $(G,\alpha)$ is a contractible group then $(G,*)$ is a quandle, with: 
$$x * y \,  =  \, x \alpha(x^{-1} y)$$

Contractible groups are particular examples of conical groups. In
\cite{buligairq} we proved that conical groups, as well as some symmetric
spaces,  can be described as emergent algebras, coming from uniform idempotent
right quasigroups.

\section{Motivation: emergent algebras}

A  differentiable algebra, is an algebra (set of operations $\mathcal{A}$) over 
a manifold $X$ with the property that all the operations 
of the algebra are differentiable with respect to the manifold structure 
of $X$. Let us denote by $\mathcal{D}$ the differential structure of the 
manifold $X$. 

From a more computational viewpoint, we may think about the calculus which can
be done in  a differentiable algebra as being generated by the elements of 
a toolbox with two compartments:  
\begin{enumerate}
\item[-]  $\mathcal{A}$ contains the algebraic
information, that is the operations of the algebra, as well as  algebraic 
relations (like for example "the operation $*$ is associative", or 
"the operation $*$ is commutative", and so on), 
\item[-] $\mathcal{D}$ contains the differential structure informations, 
that is the information needed in order to 
formulate the statement "the function $f$ is differentiable", 
\item[-]  the compartments $\mathcal{A}$ and  $\mathcal{D}$ are compatible, 
in the  that any operation from $\mathcal{A}$ is differentiable 
according to $\mathcal{D}$. 
\end{enumerate}

In the paper \cite{buligairq} we proposed the notion of a emergent algebra as a
generalization of a differentiable algebra. Computations in a 
 emergent algebra (short name for
 a uniform idempotent right quasigroup, definition \ref{deftop}) are generated
 by a class $\mathcal{E}$ of operations and relations from which a algebra 
 $\mathcal{A}$ and a generalization of a differentiable structure $\mathcal{D}$ 
 "emerge". The meaning of this emergence is the following: all elements of 
 $\mathcal{A}$ and $\mathcal{D}$ (algebraic operations, relations, differential 
 operators, ...) are constructed by finite or virtually infinite "recipes", 
 which can be implemented by some class of circuits made by very simple 
 gates (the operations in the uniform idempotent right quasigroup). 
An emergent algebra space is then described by:
\begin{enumerate}
\item[-] a class of transistor-like gates (that is binary operations), 
with in/out ports labeled by points of the space and a 
internal state variable which can be interpreted as "scale". 
\item[-]  a class of elementary circuits made of such gates 
(these are the "generators" of the emergent algebra). The elementary circuits are in fact 
certain ternary operations constructed from the operations in the uniform
idempotent right quasigroup. They have the property that the output 
converges as the scale goes to zero, uniformly with respect to the input.     
\item[-]  a class of equivalence rules saying that some simple assemblies of 
elementary circuits have equivalent function (these are the "relations" of the 
emergent algebra).  
\end{enumerate}

We shall explore in more detail this point of view, concentrating on 
braided sets which are also spaces with dilations (dilatation structures).

\section{$\Gamma$-idempotent right quasigroups}
\label{sec3}

\begin{definition} 
We use the operations of a irq to define the sum, difference and inverse operations of the irq:  for 
any $x,u,v \in X$ 
\begin{enumerate}
\item[(a)] the difference operation is $\displaystyle (xuv) \, = \, \left( x \,
\circ \, u \right) \, \bullet \, \left( x \, \circ \, v \right)$. 
 By fixing the first variable $x$ we obtain the difference operation based at $x$: 
$\displaystyle v \, -^{x} \, u \, =  \, dif^{x}(u,v) \, = \, (xuv)$. 
\item[(b)] the sum operation is 
$\displaystyle )xuv( \, = \, x \, \bullet \left( \left( x \, \circ \, u \right) \, \circ \, 
 v \right)$. 
  By fixing the first variable $x$ we obtain the sum operation based at $x$: 
$\displaystyle u \, +^{x} \, v \, =  \, sum^{x}(u,v) \, = \, )xuv($.  
\item[(a)] the inverse operation is  
$\displaystyle inv(x,u) \, = \, \left( x \, \circ \, u \right) \, \bullet \,  x
$. By fixing the first variable $x$ we obtain the inverse  operator based 
at $x$: $\displaystyle  -^{x} \, u \, =  \, inv^{x} u \, = \, inv(x,u)$. 
\end{enumerate}
 For any $\displaystyle 
k \in \mathbb{Z}^{*} = \mathbb{Z} \setminus \left\{ 0 \right\}$  we define also 
the following operations: 
\begin{enumerate}
\item[-] $\displaystyle x \, \circ_{1}\,  u \, = \, x \, \circ \, u$,  $\displaystyle x \, \bullet_{1}\,  u \, = \, x \, \bullet \, u$, 
\item[-]  for any $k > 0$ let 
$\displaystyle  x \, \circ_{k+1}\,  u \, = \, x \, \circ \left(x \circ_{k} \, u \right)$ and 
$\displaystyle  x \, \bullet_{k+1}\,  u \, = \, x \, \bullet \left(x \bullet_{k} \, u \right)$, 
\item[-] for any $k < 0$ let 
$\displaystyle  x \, \circ_{k}\,  u \, = \, x \bullet_{-k} \, u$ and 
$\displaystyle  x \, \bullet_{k}\,  u \, = \, x \circ_{-k} \, u $. 
\end{enumerate}

\label{dplay}
\end{definition}

For any $\displaystyle k \in \mathbb{Z}^{*}$  the triple 
$\displaystyle (X, \circ_{k} , \bullet_{k} )$ is a irq. We denote  the difference, sum and inverse operations of $\displaystyle (X, \circ_{k} , \bullet_{k} )$  by the same symbols as 
the ones used for  $(X, \circ , \bullet )$,  with  a subscript "$k$".

In \cite{buligairq} we introduced idempotent right quasigroups and then iterates
of the operations indexed by a parameter $\displaystyle k \in \mathbb{N}$. This
was done in order to simplify the notations mostly. Here, in the presence of the
group $\Gamma$, we might define a $\Gamma$-irq.

\begin{definition}
Let $\Gamma$ be a commutative group. A $\Gamma$-idempotent right quasigroup 
is a set $X$ with a function $\displaystyle \varepsilon \in \Gamma \mapsto 
\circ_{\varepsilon}$ such that $\displaystyle (X, \circ_{\varepsilon})$ is a irq
and moreover for any $\varepsilon, \mu \in \Gamma$ and any $x, y \in X$ we have 
$$x \, \circ_{\varepsilon} \, \left( x \, \circ_{\mu} \, y \right) \, = \, 
x \, \circ_{\varepsilon \mu} \, y$$
\label{defgammairq}
\end{definition}

It is then obvious that if $(X, \circ)$ is a irq then $\displaystyle 
(X, k \in \mathbb{Z} \mapsto \circ_{k})$ is a $\mathbb{Z}$-irq (we define 
$\displaystyle x \, \circ_{0} \, y \, = \, y$).

The following is a slight modification of proposition 3.4 and point (k)
proposition 3.5 \cite{buligairq}, for
the case of $\Gamma$-irqs (the proof of this proposition is almost 
identical, with obvious modifications, with the proof of the original
proposition).

\begin{proposition}
Let $\displaystyle (X,\circ_{\varepsilon})_{\varepsilon \in \Gamma}$ be
a $\Gamma$-irq. Then we have the relations: 
\begin{enumerate}
\item[(a)] $\displaystyle \left( u \, +^{x}_{\varepsilon} \, v \right) \,
-^{x}_{\varepsilon} \, u \, = \, v $
\item[(b)] $\displaystyle u \, +^{x}_{\varepsilon} \, \left( v \, -^{x}_{\varepsilon} \, u \right) \, = \, v $
\item[(c)] $\displaystyle v \, -^{x}_{\varepsilon} \, u \, = \, \left(-^{x}_{\varepsilon} u\right)  \, +^{x \circ u}_{\varepsilon} \, v  $
\item[(d)] $\displaystyle -^{x\circ u}_{\varepsilon} \, \left( -^{x}_{\varepsilon} \, u \right) \, = \, u $
\item[(e)] $\displaystyle u \, +^{x}_{\varepsilon} \, \left( v \, +^{x\circ u}_{\varepsilon} \, w \right) \, = \, \left( u \, +^{x}_{\varepsilon} \, v \right) \, +^{x}_{\varepsilon} \, w $
\item[(f)] $\displaystyle  -^{x}_{\varepsilon} \, u \, = \,  x \, -^{x}_{\varepsilon} \, u $
\item[(g)] $\displaystyle  x \, +^{x}_{\varepsilon} \, u \, = \,  u $
\item[(k)] for any $\displaystyle \varepsilon, \mu \in \mathbb{Z}^{*}$ and any $x, u , v \in X$ we have the distributivity property: 
$$\displaystyle (x \circ_{\mu} v) \, -_{\varepsilon}^{x} \, ( x \circ_{\mu} u)
\, = \, \left( x \circ_{\varepsilon \mu} u \right) \, 
\circ_{\mu} \, \left( v \, -_{\varepsilon \mu}^{x} \, u \right)$$
\end{enumerate}
\label{pplay}
\end{proposition}

\section{Uniform idempotent right quasigroups}

Let $\Gamma$ be a topological commutative group. We suppose that 
$\Gamma$ as a topological space is separable.

\begin{definition}
Let $(X,\tau)$ be a topological space. $\tau$ is the collection of open sets in
$X$. A filter in $(X,\tau)$ is a function $\mu: \tau \rightarrow \left\{0,1\right\}$ such that: 
\begin{enumerate}
\item[(a)] $\mu(X) \ = \ 1$, 
\item[(b)] for any $A,B \in \tau$, if $A \subset B$ then $\mu(A) \leq \mu(B)$,
\item[(c)] for any $A,B \in \tau$ we have $\displaystyle \mu(A\cup B) + \mu(A\cap B) \geq \mu(A) + \mu(B)$.
\end{enumerate}
An absolute of  a separable topological commutative group $\Gamma$ is a
class $Abs(\Gamma)$  of filters $\mu$ in $\Gamma$ with the properties: 
\begin{enumerate}
\item[(i)] for any $\varepsilon \in \Gamma$ there are $A \in \tau$ and $\mu 
\in Abs(\Gamma)$ such that $\mu(A) = 1$ and $x \not \in A$, 
\item[(ii)] for any $\mu, \mu' \in Abs(\Gamma)$ there is $A \in \tau$ such that 
$\mu(A) > \mu'(A)$, 
\item[(iii)] for any $\varepsilon \in \Gamma$ and $\mu \in Abs(\Gamma)$ the 
transport of $\mu$ by $\varepsilon$, defined as $\varepsilon \, \mu (A) = 
\mu(\varepsilon A)$, belongs to $Abs(\Gamma)$. 
\end{enumerate}
Let $f: \Gamma \rightarrow (X,\tau)$ be a function from $\Gamma$ to a
separable topological space, let  $Abs(\Gamma)$ be an absolute of $\Gamma$, and 
$\mu \in Abs(\Gamma)$. We say that $f$ converges to $z \in X$ as 
$\varepsilon$ goes to $\mu$ if for any open set $A$ in $X$ with $z \in A$ we have
$\displaystyle \mu(f^{-1}(A)) = 1$. We write: 
$$\lim_{\varepsilon \rightarrow \mu} f(\varepsilon) \, = \, z$$
\label{dfilta}
\end{definition}

For example, if $\Gamma = (0,+\infty)$ with multiplication, then 
$\displaystyle Abs(\Gamma) = \left\{ 0 \right\}$ is an absolute, where "$0$" is the filter
defined by $0(A) = 1$ if and only if the number $0$ belongs to the closure of 
$A$ in $\mathbb{R}$. Also, $Abs(\Gamma) = \left\{ 0, \infty \right\}$ is an
absolute, where "$\infty$" is the filter defined by: $\infty(A) = 1$ if and only
if $A$ is unbounded. 

Let $\Gamma$ be a commutative separable topological group,  
 $\chi: \Gamma \rightarrow (0,+\infty)$ a continuous morphism and $Abs((0,+\infty))$ an 
 absolute of $(0,+\infty)$. Let  $Abs(\Gamma)$ be the class of filters on 
 $\Gamma$ constructed like this: $\mu \in Abs(\Gamma)$ if there exists 
 $\alpha \in Abs((0,+\infty))$ such that for any open set $A$ in $\Gamma$, 
 $\mu(A) = 1$ if there is an open set $B \subset (0,+ \infty)$  with 
 $\displaystyle \chi^{-1}(B) \subset A$ and $\alpha(B) = 1$. Then $Abs(\Gamma)$
 is an absolute of $\Gamma$. 
 
 Another example: let $\displaystyle \Gamma_{0}$ be a topological separable 
 commutative group,let $G$ be a finite commutative group and let 
 $\displaystyle \Gamma = \Gamma_{0} \times G$. We think now about $G$ and 
 $\displaystyle \Gamma_{0}$ as being 
 subgroups of $\Gamma$.  Let $\displaystyle
 Abs(\Gamma_{0})$ be an absolute of $\displaystyle \Gamma_{0}$. We construct 
 $Abs(\Gamma)$ as the collection of all filters $\mu$ on $\Gamma$ such that 
 there is $g \in G$ with $\displaystyle g \mu \in Abs(\Gamma_{0})$. Then
 $Abs(\Gamma)$ is an absolute of $\Gamma$.

\begin{definition}
A $\Gamma$-uniform irq $(X, *, \backslash)$ is a separable uniform  
space $X$ endowed with continuous irq operations $*$, $\backslash$ such that: 
\begin{enumerate}
\item[(C)] the operation $*$ is compactly contractive: for each compact set 
$K \subset X$ and open set $U \subset X$, with $x \in U$, there is an open set 
$\displaystyle A(K,U) \subset \Gamma$ with $\mu(A) = 1$ for any $\mu \in 
Abs(\Gamma)$ and for any $u \in K$ and
$\varepsilon \in A(K,U)$, we have $\displaystyle x *_{\varepsilon} u \in U$; 
\item[(D)] the following limits exist for any $\mu \in Abs(\Gamma)$ 
$$ \lim_{\varepsilon \rightarrow \mu} v \, -_{\varepsilon}^{x} \, u \, = \, 
v \, -_{\infty}^{x} \, u \quad , \quad \lim_{\varepsilon \rightarrow \mu} u
\,+_{\varepsilon}^{x}
\, v \,  = \, u \,+_{\infty}^{x} \, v $$
and are uniform with respect to $x, u, v$ in a compact set. 
\end{enumerate}
\label{deftop}
\end{definition}

The main property of a uniform irq is the following. It is a consequence of 
relations from proposition \ref{pplay}. 

\begin{theorem}
Let $(X, *, \backslash)$ be a uniform irq. Then for any $x \in X$  the operation 
$\displaystyle (u,v) \mapsto u \,+_{\infty}^{x} \, v$ gives $X$ the structure of
a conical group with the dilatation $\displaystyle u \mapsto x *_{\varepsilon} u$.
\label{mainthm}
\end{theorem}

Conical groups are described in the next section.

\section{Conical groups are distributive uniform irqs}
\label{secon}

For a dilatation structure (see section \ref{dilats}) the metric tangent spaces   have a group structure 
which is compatible with dilatations. This structure, of a  group with dilatations, is interesting 
by itself. The notion has been introduced in \cite{buliga2}, \cite{buligadil1}; we describe  it 
further.

Let $\Gamma$ be a topological commutative groups with an absolute $Abs(\Gamma)$.

\begin{definition}
A group with dilatations $(G,\delta)$ is a  topological group $G$  with  an action 
of $\Gamma$ (denoted by $\delta$), on $G$ such that for any $\mu \in Abs(\Gamma)$
\begin{enumerate}
\item[H0.] the limit  $\displaystyle \lim_{\varepsilon \rightarrow \mu} 
\delta_{\varepsilon} x  =  e$ exists and is uniform with respect to $x$ in a compact neighbourhood of the identity $e$.
\item[H1.] the limit
$$\beta (x,y)  =  \lim_{\varepsilon \rightarrow \mu} \delta_{\varepsilon}^{-1}
\left((\delta_{\varepsilon}x) (\delta_{\varepsilon}y ) \right)$$
is well defined in a compact neighbourhood of $e$ and the limit is uniform.
\item[H2.] the following relation holds
$$ \lim_{\varepsilon \rightarrow \mu} \delta_{\varepsilon}^{-1}
\left( ( \delta_{\varepsilon}x)^{-1}\right)  =  x^{-1}$$
where the limit from the left hand side exists in a neighbourhood of $e$ and is uniform with respect to $x$.
\end{enumerate}
\label{defgwd}
\end{definition}

\begin{definition}
A conical group $(N, \delta)$ is a   group with dilatations  such that for any $\varepsilon \in \Gamma$  the dilatation 
 $\delta_{\varepsilon}$ is a group morphism. 
\end{definition}

A conical group is the infinitesimal version of a group with 
dilatations (\cite{buligadil1} proposition 2).

\begin{proposition}
Under the hypotheses H0, H1, H2,  $\displaystyle (G,\beta, \delta)$,  is a conical group, with operation 
$\displaystyle \beta$,  dilatations $\delta$.
\label{here3.4}
\end{proposition}

One particular case is the one of contractible groups, definition \ref{defunu}, 
which are also normed groups. Indeed, in this case we may take 
$\Gamma = \mathbb{Z}$.

Locally compact conical groups are  locally compact groups admitting 
a contractive automorphism group. We begin with  the
definition of a contracting automorphism group \cite{siebert}, definition 5.1. 

\begin{definition}
Let $G$ be a locally compact group. An automorphism group on $G$ is a family 
$\displaystyle T= \left( \tau_{t}\right)_{t >0}$ in $Aut(G)$, such that 
$\displaystyle \tau_{t} \, \tau_{s} = \tau_{ts}$ for all $t,s > 0$. 

The contraction group of $T$ is defined by 
$$C(T) \ = \ \left\{ x \in G \mbox{ : } \lim_{t \rightarrow 0} \tau_{t}(x) = e
\right\} \quad .$$
The automorphism group $T$ is contractive if $C(T) = G$. 
\end{definition}

Next is proposition 5.4 \cite{siebert}, which gives a description of locally
compact groups which admit a contractive automorphism group. 

\begin{proposition}
For a locally compact group $G$ the following assertions are equivalent: 
\begin{enumerate}
\item[(i)] $G$ admits a contractive automorphism group;
\item[(ii)] $G$ is a simply connected Lie group whose Lie algebra admits a 
positive graduation.
\end{enumerate}
\label{psiebert}
\end{proposition}

The proof of the next proposition is an easy  application of 
the previously explained facts. 

\begin{proposition}
Let $(G,\delta)$ be a locally compact conical group. Then the associate 
irq  $(G,*)$ is an uniform irq. 
\label{pfirstop}
\end{proposition}

A particular class of locally compact groups which admit a contractive
automorphism group is  made by Carnot groups. They are related to sub-riemannian or 
Carnot-Carath\'eodory geometry, which is the study of non-holonomic manifolds
endowed with a Carnot-Carath\'eodory distance. Non-holonomic spaces were
discovered in 1926 by  G. Vr\u anceanu \cite{vra1},
\cite{vra2}. The  Carnot-Carath\'eodory distance on a non-holonomic space is 
inspired by Carath\'eodory \cite{cara} work from 1909   on the mathematical formulation of
thermodynamics. Such spaces appear in applications
to thermodynamics, to the mechanics of non-holonomic systems, in the study of
hypo-elliptic operators cf. H\"ormander \cite{hormander}, in harmonic analysis on
homogeneous cones cf. Folland, Stein \cite{fostein}, and
as boundaries of CR-manifolds.

The following result is a slight modification of  \cite{buligairq}, theorem 6.1,
consisting in the replacement of "contractible" by "conical" in the statement of
the theorem. 

\begin{theorem}
Let $(G,\alpha)$ be a locally compact conical  group and $G(\alpha)$ be the 
associated uniform irq. Then the irq is distributive, namely it satisfies the
relation: for any $\varepsilon, \lambda \in \Gamma$ 
\begin{equation}
x *_{\varepsilon} \left( y *_{\lambda} z \right) \, = \, \left( x
*_{\varepsilon} y \right) *_{\lambda} \left( 
x *_{\varepsilon} z \right) 
\label{distributive}
\end{equation}

Conversely, let $(G, *)$ be a distributive uniform irq. Then there 
is a group operation on $G$ (denoted multiplicatively), with neutral element
$e$, such that:
\begin{enumerate}
\item[(i)] $\displaystyle xy \, = \, x +^{e}_{\infty} y$ for any $x, y \in G$, 
\item[(ii)] for any $x, y , z \in G$ we have $\displaystyle (xyz)_{\infty} = 
x y^{-1} z$, 
\item[(iii)] for any $x, y \in G$ we have $\displaystyle 
x *_{\varepsilon} y \, = \, x (e *_{\varepsilon} (x^{-1} y))$. 
\end{enumerate}
In conclusion  there is a bijection between distributive $\Gamma$-uniform irqs and 
conical groups. 
\label{pgroudlin}
\end{theorem}

\section{Normed uniform irqs are dilatation structures} 
\label{dilats}

For simplicity we shall list the   axioms of  a dilatation structure $(X,d,\delta)$ without concerning about
domains and codomains of dilatations. For the full definition of dilatation
structures, as well as for their main properties and examples, see
\cite{buligadil1}, \cite{buligadil2}, \cite{buligasr}. The notion appeared from
my efforts to understand  the last section of the paper  
\cite{bell} (see also \cite{pansu}, \cite{gromovsr}, \cite{marmos1},
\cite{marmos2}).

Let $\Gamma$ be a topological commutative groups with an absolute $Abs(\Gamma)$ 
and with a morphism $\mid \cdot \mid : \Gamma \rightarrow (0,+\infty)$ such that
for any $\mu \in Abs(\Gamma)$ 
$$\lim_{\varepsilon \rightarrow \mu} \mid \varepsilon \mid \, = \, 0$$

\begin{definition}
A triple $(X, d, \delta)$ is a dilatation structure if $(X, d)$ is a locally 
compact metric space and the dilatation field 
 $$\delta: \Gamma \times \left\{ (x,y) \in X \times X \mbox{ : } y \in
 dom(\varepsilon, x) \right\} \rightarrow X \quad , \quad \delta(\varepsilon, x,
 y) \, = \, \delta^{x}_{\varepsilon} y $$
 gives to $X$ the structure of a uniform idempotent right quasigroup over 
 $\Gamma$ (definition
 \ref{deftop}), with the operation: for any $\varepsilon \in \Gamma$  
 $$x \, *_{\varepsilon} \, y \, = \, \delta_{\varepsilon}^{x} y$$
 Moreover, the distance is compatible with the dilatations, in the sense: 
 
\begin{enumerate}
\item[A1.]  the uniformity on $(X,\delta)$ is the one induced by the
distance $d$,

\item[A2.] There is $A > 1$ such that  for any 
$x$ there exists
 a  function $\displaystyle (u,v) \mapsto d^{x}(u,v)$, defined for any
$u,v$ in the closed ball (in distance d) $\displaystyle
\bar{B}(x,A)$, such that for any $\mu \in Abs(\Gamma)$
$$\lim_{\varepsilon \rightarrow \mu} \quad \sup  \left\{  \mid \frac{1}{\mid 
\varepsilon \mid} d(\delta^{x}_{\varepsilon} u,
\delta^{x}_{\varepsilon} v) \ - \ d^{x}(u,v) \mid \mbox{ :  } u,v \in \bar{B}_{d}(x,A)\right\} \ =  \ 0$$
uniformly with respect to $x$ in compact set. 
Moreover the uniformity induced by $d^{x}$ is the same as the uniformity induced
by $d$, in particular $\displaystyle d^{x}(u,v) = 0$ implies $u = v$.

\end{enumerate}

\label{defweakstrong}
\end{definition}

In order to make connection with the original definition of a dilatation
structure introduced and studied in \cite{buligadil1}, \cite{buligadil2}, we 
shall relate the notations used here with the original ones. 
The operations induced by the uniform irq structure on $X$ are: 
$$ v \, -_{\varepsilon}^{x} \, u \, = \, \Delta^{x}_{\varepsilon}(u,v) \, = \, 
\delta_{\varepsilon^{-1}}^{\delta^{x}_{\varepsilon}u} \delta^{x}_{\varepsilon} v$$ 
$$ u \,+_{\varepsilon}^{x}
\, v \,  = \, \Sigma^{x}_{\varepsilon}(u,v) \, = \, 
\delta_{\varepsilon^{-1}}^{x} \delta_{\varepsilon}^{\delta^{x}_{\varepsilon}u} 
v $$
where $\displaystyle \Delta^{x}_{\varepsilon}$, $\Sigma^{x}_{\varepsilon}$ 
are the approximate difference, respectively approximate sum operations induced 
by dilatation structures. Similarily we have the following correspondence of
notations: 
$$ v \, -_{\infty}^{x} \, u \, = \, \Delta^{x}(u,v) \, = \, \lim_{\varepsilon
\rightarrow  \mu} 
\delta_{\varepsilon^{-1}}^{\delta^{x}_{\varepsilon}u} \delta^{x}_{\varepsilon} v$$ 
$$ u \,+_{\varepsilon}^{x}
\, v \,  = \, \Sigma^{x}(u,v) \, = \, \lim_{\varepsilon
\rightarrow  \mu} 
\delta_{\varepsilon^{-1}}^{x} \delta_{\varepsilon}^{\delta^{x}_{\varepsilon}u} 
v $$

The conclusion is therefore that adding a distance in the story of uniform 
irqs gives us the notion of a dilatation structure. 

We go a bit into details. 

\begin{proposition}
Let $(X,d, \delta)$ be a dilatation structure,  $x \in X$, and let 
$$\delta^{x}_{\varepsilon} \, d (u, v) \, = \, \frac{1}{\mid \varepsilon \mid}
\, d( \delta^{x}_{\varepsilon} u , \delta^{x}_{\varepsilon} v )$$
 Then the net of metric spaces $\displaystyle (\bar{B}_{d}(x,A),
 \delta^{x}_{\varepsilon} d)$ converges in the Gromov-Hausdorff sense to the 
 metric space $\displaystyle (\bar{B}_{d}(x,A), d^{x})$. Moreover this metric
 space is a metric cone, in the following sense: for any $\lambda \in \Gamma$ 
  we have
 $$d^{x} ( \delta^{x}_{\lambda} u , \delta^{x}_{\lambda} v ) \, = \, 
 \mid \lambda \mid 
 \, d^{x}(u,v)$$
 \label{pizo}
 \end{proposition}
 
 \paragraph{Proof.} 
 The first part of the proposition is just a reformulation  of axiom A2, 
 without the condition of uniform convergence. For the second part 
 remark  that 
  
 $$\frac{1}{\mid \varepsilon \mid} d( \delta^{x}_{\varepsilon} \, 
 \delta^{x}_{\lambda} u , \delta^{x}_{\varepsilon} \, 
 \delta^{x}_{\lambda} v ) \, = \, \mid \lambda \mid \, \frac{1}{\mid 
 \varepsilon \lambda \mid} d( \delta^{x}_{\varepsilon \lambda}  u , 
 \delta^{x}_{\varepsilon \lambda}  v ) $$
 Therefore if we pass to the limit 
 with $\varepsilon \rightarrow \mu$ in these two relations we get the desired
 conclusion. \quad $\square$

Particular examples of dilatation structures are given by normed groups with
dilatations. 

\begin{definition} A normed group with dilatations $(G, \delta, \| \cdot \|)$ is a 
group with dilatations  $(G, \delta)$ endowed with a continuous norm  
function $\displaystyle \|\cdot \| : G \rightarrow \mathbb{R}$ which satisfies 
(locally, in a neighbourhood  of the neutral element $e$) the properties: 
 \begin{enumerate}
 \item[(a)] for any $x$ we have $\| x\| \geq 0$; if $\| x\| = 0$ then $x=e$, 
 \item[(b)] for any $x,y$ we have $\|xy\| \leq \|x\| + \|y\|$, 
 \item[(c)] for any $x$ we have $\displaystyle \| x^{-1}\| = \|x\|$, 
 \item[(d)] the limit 
$\displaystyle \lim_{\varepsilon \rightarrow \mu} \frac{1}{\mid\varepsilon \mid} \| \delta_{\varepsilon} x \| = \| x\|^{N}$ 
 exists, is uniform with respect to $x$ in compact set, 
 \item[(e)] if $\displaystyle \| x\|^{N} = 0$ then $x=e$.
  \end{enumerate}
  \label{dnco}
  \end{definition}

In a normed group with dilatations we have a natural left invariant distance given by
\begin{equation}
d(x,y) = \| x^{-1}y\| \quad . 
\label{dnormed}
\end{equation}
Any normed group with dilatations has an associated dilatation structure on it.  In a group with dilatations $(G, \delta)$  we define dilatations based in any point $x \in G$ by 
 \begin{equation}
 \delta^{x}_{\varepsilon} u = x \delta_{\varepsilon} ( x^{-1}u)  . 
 \label{dilat}
 \end{equation}

The following result is theorem 15 \cite{buligadil1}. 

\begin{theorem}
Let $(G, \delta, \| \cdot \|)$ be  a locally compact  normed  group with dilatations. Then $(G, d, \delta)$ is 
a dilatation structure, where $\delta$ are the dilatations defined by (\ref{dilat}) and the distance $d$ is induced by the norm as in (\ref{dnormed}). 
\label{tgrd}
\end{theorem}

The general theorem \ref{mainthm} has a stronger conclusion in the case of
dilatation structures, namely "conical groups" are replaced by "normed conical
groups". 

\section{Differentiability}

We have seen that to any uniform irq we can associate a bundle of contractible
groups  $\displaystyle x \in X \mapsto (X, +^{x}_{\infty}, x * )$. This
bundle can be seen as a tangent bundle, namely: to any $x \in X$ is
associated a conical group with $x$ as neutral element, which 
 is the tangent space at $x$. We shall denote it by $\displaystyle T^{x} X$.

 In the particular case of a manifold, this is indeed a correct definition 
  in the following sense: if we look to a small 
 portion of the manifold then we know that there is a chart of this small 
 portion, which puts it in bijection with an open set in $\mathbb{R}^{n}$. We
 have seen that we can associate to $\mathbb{R}^{n}$ a uniform irq by using 
 as the operation $*$ a homothety with fixed ratio $\varepsilon < 1$. This 
 uniform irq is transported on the manifold by the chart. If we ignore the facts
  that we are working not with the whole manifold, but with a small part of it, 
 and not with $\mathbb{R}^{n}$, but with a open set, then indeed we may 
 identify, for any point $x$ in the manifold, a neighbourhood of the point with 
 a neighbourhood of the tangent space at the point, such that the operation 
 of addition of vectors in the tangent space at $x$ is just the operation 
 $\displaystyle +^{x}_{\infty}$ and scalar multiplication by (any integer power 
 of) $\varepsilon$ is just $\displaystyle u \mapsto  x *_{\varepsilon} u$ . 
 
 The same is true in the more complex situation of a sub-riemannian manifold, as
 shown in \cite{buligasr}, in the sense that (locally) we may associate to 
 each point $x$ a "dilatation" of ratio $\varepsilon$, which in turn gives 
 us a structure of uniform irq. In the end we get a bundle of Carnot group 
 operations which can be seen as a tangent bundle of the sub-riemannian
 manifold. (In this case we actually have more structure given by the
 Carnot-Carath\'eodory distance, which induces also a "group norm" on each 
 Carnot group.) 
 
 A uniform irq can be seen as a generalization of a 
 differential structure. For this we give a definition of 
 differentiable functions between two uniform irqs. This definition 
 corresponds to uniform differentiability in the metric case of 
 dilatation structures, definition 16 and the comments after it in 
 \cite{buligadil1}. It is a generalization of Pansu differentiability
 \cite{pansu}.
 
 \begin{definition}
 Let $(X, *, \backslash)$ and $(Y, \circ ,  \backslash \backslash )$ be two 
$\Gamma$-uniform irqs. A function  $f: X \rightarrow  Y$ is differentiable
 if there is a function $Tf: X \times X \rightarrow Y$ such that 
 $$\lim_{k \rightarrow \infty} f(x) \, \backslash \backslash_{k} \, f (x *_{k}
 u) \, = \, Tf (x, u)$$
 uniformly with respect to $x, u$ in compact sets. 
 \label{defder}
 \end{definition}
 By abstract  nonsense the application $Tf$ has nice properties, like 
 $\displaystyle Tf(x, \cdot) : (X, +^{x}) \rightarrow (Y, +^{f(x)})$ is a
 morphism of conical groups.

\section{Sub-riemannian symmetric spaces as braided 
$\displaystyle \mathbb{R} \times \mathbb{Z}_{2}$-dilatation structures}
\label{secbraid}

Sub-riemannian symmetric spaces have been introduced in \cite{stric}, 
section 9.  We shall be interested in the description 
of sub-riemannian geometry by dilatation structures, therefore we shall 
use the same notations as in the previous paper \cite{buligasr} (see also the 
relevant citations in that paper, as well as the long paper \cite{buligadil3},
where the study of sub-riemannian geometry as a length dilatation structure is 
completed).

\begin{definition}
(adaptation of \cite{stric} definition 8.1) Let $(M,D,g)$ be a regular sub-riemannian
manifold. We say that $\Psi: M \rightarrow M$ is an infinitesimal isometry 
if $\Psi$ is $\displaystyle \mathcal{C}^{1}$ and $D\Psi$ preserves the metric
$g$. An infinitesimal isometry is regular if for any $x \in M$ and any 
tangent vector $\displaystyle u \in T_{x}M$ 
$$\Psi(exp_{x}(u)) \, = \, exp_{\Psi(x)}(D\Psi(x)u)$$
\label{dinfiso}
\end{definition}

By \cite{stric} theorem 8.2., $\displaystyle \mathcal{C}^{1}$ isometries are 
regular infinitesimal isometries and, conversely, regular infinitesimal
isometries are isometries.

An equivalent description of regular infinitesimal isometries is the following: 
they are $\mathcal{C}^{1}$ Pansu differentiable isometries.

\begin{definition}
(\cite{stric} definition 9.1) A sub-riemannian symmetric space is a regular
sub-riemannian manifold $M,D,g)$ which has a transitive Lie group $G$ of regular
infinitesimal isometries acting differentiably on $M$ such that: 
\begin{enumerate}
\item[(i)] there is a point $x \in X$ such that the isotropy subgroup $K$ 
of $x$  is compact, 
\item[(ii)] $K$ contains an element $\Psi$ such that 
$\displaystyle D\Psi(x)_{|_{D_{x}}} = - \, id$ and $\Psi$ is involutive. 
\end{enumerate}
If $G$ is a group for which (i), (ii) holds then we call $G$ an admissible
isometry group for $M$.
\label{defsrsym}
\end{definition}

\begin{theorem}
(\cite{stric} theorem 9.2) If $M$ is a sub-riemannian symmetric space and $G$ 
is
an admissible isometry group, then there exists an involution $\sigma$ of $G$ 
such that $\sigma(K) \subset K$ with the following properties 
(we write $\displaystyle \mathfrak{g} = \mathfrak{g}^{+} + \mathfrak{g}^{-}$, 
where $\displaystyle \mathfrak{g}^{+}, \mathfrak{g}^{-}$ are the subspaces of 
$\displaystyle \mathfrak{g}$ on which $D\sigma$ acts as $Id, \, - Id$): 
\begin{enumerate}
\item[(a)] $\mathfrak{g}$ is generated as a Lie algebra by a subspace
$\mathfrak{p}$ and the Lie algebra  $\mathfrak{t}$ of $K$ with $\displaystyle 
\mathfrak{p} \subset \mathfrak{g}^{-}$, $\displaystyle \mathfrak{t} \subset 
\mathfrak{g}^{+}$, 
\item[(b)] there exists a positive definite quadratic form $g$ on $\mathfrak{p}$
and $ad \, K$ maps $\mathfrak{p}$ to itself and preserves $g$. Furthermore,
$\mathfrak{p}$ may be identified with $\displaystyle D_{x}$ under the
exponential map of the Lie algebra $\mathfrak{g}$, and $g$ with the
sub-riemannian metric on $\displaystyle D_{x}$. 
\end{enumerate}

Conversely, given a Lie group $G$ and an involution $\sigma$ such that (a) and
(b) hold, then $G/K$ forms a sub-riemannian symmetric space, where
$\displaystyle D_{x_{0}} = \, exp \, \mathfrak{p}$ for the point $\displaystyle 
x_{0}$ identified with the coset $K$, and the sub-riemannian metric on 
$\displaystyle D_{x_{0}}$ is given by $g$. The bundle $D$ and its metric is then
uniquely determined by the requirement that elements of $G$ be infinitesimal 
isometries. 
\label{stricthm}
\end{theorem}

As a consequence of this theorem we see that we may endow a sub-riemannian
symmetric space,  with admissible isometry group $G$, with a (reflexive space) 
operation 
$$(x, y) \in M^{2} \mapsto \Psi(x,y) = \Psi^{x} y$$ 
such that $\Psi$ is distributive, for any $x \in X$ the map $\displaystyle 
\Psi^{x}$ satisfies (ii) definition \ref{defsrsym}, 
 and for any $g \in G$ and any $x, y \in X$ we have
$$g \left( \Psi^{x} y \right) \, = \, \Psi^{g(x)}g(y)$$

We explained in \cite{buligasr} that we can construct a dilatation structure 
over a regular sub-riemannian manifold by using adapted frames. 

Let us consider now  dilatations structures with 
$\Gamma$ isomorphic with $\displaystyle \mathbb{R} \times \mathbb{Z}_{2}$. That
means $\Gamma$ is the commutative group made by two copies of 
$(0,+ \infty)$, generated by $(0,+\infty)$ and an element $\sigma \not \in 
(0, + \infty)$, with the properties: for any $\varepsilon \in (0,+\infty)$ 
we have $\varepsilon \sigma = \sigma \varepsilon$ and $\sigma \sigma = 1$. The
absolute we take has two elements, one corresponding to $\varepsilon \rightarrow
0$ (we denote it by "$0$") and the other one is the transport by $\sigma$ of 
$0$, denoted by "$0\sigma$". The morphism $\mid \cdot \mid$ is defined by 
$$\mid \varepsilon \mid \, = \, \mid \sigma \varepsilon \mid \, = \,
\varepsilon$$

Let $(X,d, \delta)$ be a dilatation structure with respect to the group 
$\Gamma$ , absolute $Abs(\Gamma)$ and morphism $\mid \cdot \mid$ described
previously. Then for any $\varepsilon \in (0,+\infty)$ and any $x \in X$ we have
the relations:
$$\delta_{\sigma}^{x} \, \delta_{\varepsilon}^{x} \, = \,
\delta_{\varepsilon}^{x} \, \delta^{x}_{\sigma} \quad , \quad
\delta^{x}_{\sigma} \, \delta^{x}_{\sigma} \, = \, id$$

\begin{proposition}
Denote by $\displaystyle \sigma^{x} y \, = \, \delta^{x}_{\sigma} y$ and 
suppose that for any $x \in X$ the map $\displaystyle \sigma^{x}$ is not
the identity map. Then $\displaystyle \sigma^{x}$ is involutive, a isometry 
of $\displaystyle d^{x}$ and an isomorphism of the conical group $\displaystyle 
T^{x}X$.
\label{psig}
\end{proposition}

\paragraph{Proof.}
For any $x \in X$ clearly $\displaystyle \sigma^{x}$ is involutive, commutes
with dilatations $\displaystyle \delta^{x}_{\varepsilon}$ and as a consequence 
of proposition \ref{pizo}, is an isometry of of $\displaystyle d^{x}$. We need
to show that it preserves the operation $\displaystyle +^{x}_{\infty}$. We shall
work with the notations from dilatation structures. We have then, for any 
$\varepsilon \in (0,+\infty)$: 
$$\sigma^{\delta^{x}_{\varepsilon \sigma} u} \, \Delta^{x}_{\sigma \varepsilon} (u,v) \, = \, \delta^{\delta^{x}_{\sigma
\varepsilon} u}_{
\varepsilon^{-1}} \, \delta^{x}_{\sigma \varepsilon} v \, = \, 
 \Delta^{x}_{\varepsilon}(\sigma^{x} u
, \sigma^{x} v)$$
We pass to the limit with $\varepsilon \rightarrow 0$ and we get the relation: 
$$\sigma^{x} \Delta^{x}(u,v) \, = \, \Delta^{x}(\sigma^{x} u , \sigma^{x} v)$$
which shows that $\displaystyle \sigma^{x}$ is an isomorphism of $\displaystyle 
T^{x}X$. \quad $\square$

This proposition motivates us to introduce braided $\displaystyle \mathbb{R}
\times \mathbb{Z}_{2}$-dilatation structures. 

\begin{definition}
Let  $(X, d, \delta)$ be a dilatation structure, with respect to the group 
$\Gamma$ , absolute $Abs(\Gamma)$ and morphism $\mid \cdot \mid$ described
previously, and such that for any $x \in X$ the map $\displaystyle \sigma^{x}$ is not
the identity map. This dilatation structure is braided if the map 
$$(x,y) \in X^{2} \mapsto (\sigma^{x}y, x)$$
is a braided map. 
\end{definition}

\begin{theorem}
A sub-riemannian symmetric space $M$ with admissible isometry group $G$  can be endowed with a braided 
$\displaystyle \mathbb{R} \times \mathbb{Z}_{2}$-dilatation structure 
which is $G$-invariant, that is for any $g \in G$, for any $x,y \in M$, and for any 
$\varepsilon \in \Gamma$ we have 
 $$g \left( \delta^{x}_{\varepsilon} y \right) \, = \, \delta_{\varepsilon}^{g(x)}g(y)$$
\end{theorem}

\paragraph{Sketch of the proof.} 
 In the particular case of a sub-riemannian symmetric space we may obviously take the
adapted frames to be $G$-invariant, therefore we may construct a dilatation
structure (over the group $(0,+\infty)$ with multiplication) 
which is $G$-invariant. Because $\displaystyle \Psi^{x}$ satisfies (ii) definition \ref{defsrsym}, it 
follows $\displaystyle \Psi^{x}$ is differentiable in $x$ in the sense of
dilatation structures. We extend the dilatation structure to a braided one 
by defining for any $x \in X$ 
$$\sigma^{x} \, = \, T\Psi^{x}(x,
\cdot)$$
By $G$-invariance of both the dilatation structure and the operation $\Psi$ it
follows that 
$$T\Psi^{x}(x, \cdot) \, = \, \Psi^{x}$$ 
therefore  $\sigma^{x}$  commutes with  $\displaystyle
\delta^{x}_{\varepsilon}$, which ensures us that we well defined a braided 
dilatation structure. \quad $\square$

\end{document}